\documentstyle[12pt]{article}
\begin{document}
\newcommand{\text}[1]{\mbox{{\rm #1}}}  
\newcommand{\gd}{\delta}
\newcommand{\itms}[1]{\item[[#1]]} 
\newcommand{\nin}{\in\!\!\!\!\!/}
\newcommand{\sub}{\subset} 
\newcommand{\cntd}{\subseteq}     
\newcommand{\go}{\omega} 
\newcommand{\Pa}{P_{a^\nu,1}(U)} 
\newcommand{\fx}{f(x)}  
\newcommand{\fy}{f(y)} 
\newcommand{\gD}{\Delta}
\newcommand{\gl}{\lambda} 
\newcommand{\gL}{\Lambda} 
\newcommand{\half}{\frac{1}{2}} 
\newcommand{\sto}[1]{#1^{(1)}}
\newcommand{\stt}[1]{#1^{(2)}}
\newcommand{\Z}{\hbox{\sf Z\kern-0.720em\hbox{ Z}}}
\newcommand{\C}{C}
\newcommand{\singcolb}[2]{\left(\begin{array}{c}#1\\#2
\end{array}\right)} 
\newcommand{\ga}{\alpha}
\newcommand{\gb}{\beta} 
\newcommand{\gga}{\gamma}
\newcommand{\ul}{\underline} 
\newcommand{\ol}{\overline} 
\newcommand{\qed}{\kern 5pt\vrule height8pt width6.5pt depth2pt}
\newcommand{\Lrraro}{\Longrightarrow}
\newcommand{\Nb}{|\!\!/}
\newcommand{\NN}{{\rm I\!N}}
\newcommand{\bsl}{\backslash}     
\newcommand{\gt}{\theta}
\newcommand{\op}{\oplus}
\newcommand{\Op}{\bigoplus}          
\newcommand{\CR}{{\cal R}}
\newcommand{\tr}{\bigtriangleup}
\newcommand{\grr}{\omega_1} 
\newcommand{\ben}{\begin{enumerate}}
\newcommand{\een}{\end{enumerate}}
\newcommand{\ndiv}{\not\mid}
\newcommand{\bab}{\bowtie}
\newcommand{\hal}{\leftharpoonup}
\newcommand{\har}{\rightharpoonup}
\newcommand{\ot}{\otimes}
\newcommand{\OT}{\bigotimes}
\newcommand{\bwe}{\bigwedge}
\newcommand{\gep}{\varepsilon}
\newcommand{\gs}{\sigma} 
\newcommand{\rbraces}[1]{\left( #1 \right)}
\newcommand{\bbox}{$\;\;\rule{2mm}{2mm}$}
\newcommand{\sbraces}[1]{\left[ #1 \right]}
\newcommand{\bbraces}[1]{\left\{ #1 \right\}}
\newcommand{\OO}{_{(1)}}
\newcommand{\TT}{_{(2)}}
\newcommand{\FF}{_{(3)}}
\newcommand{\minus}{^{-1}}
\newcommand{\CV}{\cal V} 
\newcommand{\CVs}{\cal{V}_s} 
\newcommand{\un}{U_q(sl_n)'}
\newcommand{\on}{O_q(SL_n)'}
\newcommand{\slq}{U_q(sl_2)}
\newcommand{\olq}{O_q(SL_2)}
\newcommand{\UU}{U_{(N,\nu,\go)}}
\newcommand{\HH}{H_{n,q,N,\nu}} 
\newcommand{\ct}{\centerline}
\newcommand{\bs}{\bigskip}
\newcommand{\qua}{\rm quasitriangular}   
\newcommand{\ms}{\medskip}
\newcommand{\noin}{\noindent}
\newcommand{\mat}[1]{$\;{#1}\;$}
\newcommand{\raro}{\rightarrow}
\newcommand{\map}[3]{{#1}\::\:{#2}\raro{#3}}
\newcommand{\alg}{{\rm Alg}}
\def\newtheorems{\newtheorem{theorem}{Theorem}[subsection]
                 \newtheorem{cor}[theorem]{Corollary}
                 \newtheorem{prop}[theorem]{Proposition}
                 \newtheorem{lemma}[theorem]{Lemma}
                 \newtheorem{defn}[theorem]{Definition}
                 \newtheorem{Theorem}{Theorem}[section]
                 \newtheorem{Corollary}[Theorem]{Corollary}
                 \newtheorem{Proposition}[Theorem]{Proposition}
                 \newtheorem{Lemma}[Theorem]{Lemma}
                 \newtheorem{Defn}[Theorem]{Definition}
                 \newtheorem{Example}[Theorem]{Example}
                 \newtheorem{Remark}[Theorem]{Remark} 
                 \newtheorem{claim}[theorem]{Claim}
                 \newtheorem{sublemma}[theorem]{Sublemma}
                 \newtheorem{example}[theorem]{Example}
                 \newtheorem{remark}[theorem]{Remark}
                 \newtheorem{question}[theorem]{Question}
                 \newtheorem{conjecture}{Conjecture}[subsection]}
\newtheorems
\newcommand{\proof}{\par\noindent{\bf Proof:}\quad}
% Matrix ():
\newcommand{\dmatr}[2]{\left(\begin{array}{c}{#1}\\
                            {#2}\end{array}\right)}
\newcommand{\doubcolb}[4]{\left(\begin{array}{cc}#1&#2\\
#3&#4\end{array}\right)}
\newcommand{\qmatrl}[4]{\left(\begin{array}{ll}{#1}&{#2}\\
                            {#3}&{#4}\end{array}\right)}
\newcommand{\qmatrc}[4]{\left(\begin{array}{cc}{#1}&{#2}\\
                            {#3}&{#4}\end{array}\right)}
\newcommand{\qmatrr}[4]{\left(\begin{array}{rr}{#1}&{#2}\\
                            {#3}&{#4}\end{array}\right)}
\newcommand{\smatr}[2]{\left(\begin{array}{c}{#1}\\
                            \vdots\\{#2}\end{array}\right)}

% Matrix []:
\newcommand{\ddet}[2]{\left[\begin{array}{c}{#1}\\
                           {#2}\end{array}\right]}
\newcommand{\qdetl}[4]{\left[\begin{array}{ll}{#1}&{#2}\\
                           {#3}&{#4}\end{array}\right]}
\newcommand{\qdetc}[4]{\left[\begin{array}{cc}{#1}&{#2}\\
                           {#3}&{#4}\end{array}\right]}
\newcommand{\qdetr}[4]{\left[\begin{array}{rr}{#1}&{#2}\\
                           {#3}&{#4}\end{array}\right]}

% Matrix {}:
\newcommand{\qbracl}[4]{\left\{\begin{array}{ll}{#1}&{#2}\\
                           {#3}&{#4}\end{array}\right.}
\newcommand{\qbracr}[4]{\left.\begin{array}{ll}{#1}&{#2}\\
                           {#3}&{#4}\end{array}\right\}}

%Introduction
\title{Some Properties of Finite-Dimensional 
Semisimple Hopf Algebras}
\author{Pavel Etingof and Shlomo Gelaki
\\Department of Mathematics\\
Harvard University\\Cambridge, MA 02138}
\date{November 26, 1997}
\maketitle
Kaplansky conjectured that if $H$ is a finite-dimensional semisimple
Hopf algebra over an algebraically closed field $k$ of characteristic $0$,
then $H$ is of Frobenius type (i.e.
if $V$ is an irreducible representation of $H$ then $dimV$ divides
$dimH$) [Ka]. It was proved that the conjecture is true for
$H$ of dimension $p^n,$ $p$ prime \cite{mw}, and that if $H$
has a $2-$dimensional representation then $dimH$ is even \cite{nr}. 

In this paper we first prove in {\bf Theorem \ref{main2}} that if $V$ is
an irreducible representation of
$D(H),$ the Drinfeld double of any finite-dimensional semisimple Hopf
algebra $H$ over $k,$ then $dimV$ divides $dimH$ (not just 
$dimD(H)=(dimH)^2$). In doing this we use the theory of modular tensor 
categories (in particular Verlinde formula). We then 
use it to prove in {\bf Theorem \ref{main3}} that Kaplansky's conjecture is 
true for 
finite-dimensional semisimple quasitriangular Hopf algebras over $k.$
As a result we prove easily in {\bf Theorem \ref{p}} that 
Kaplansky's conjecture \cite{ka} on prime dimensional 
Hopf algebras over $k$ is true by passing to their Drinfeld doubles
(compare with \cite{z}). 

Second, we use a theorem of Deligne \cite{de} to prove in {\bf Theorem 
\ref{main1}} that triangular 
semisimple Hopf algebras over 
$k$ are equivalent to group algebras as quasi-Hopf algebras \cite{dr2}.

\section{Quasitriangular Semisimple Hopf Algebras are of\\ Frobenius Type}
Throughout this paper, unless otherwise is specified, $k$ will denote an
algebraically closed field of characteristic $0.$ 

Let $(H,R)$ be a finite-dimensional quasitriangular Hopf algebra over
$k,$ and write $R=\sum_i a_i\ot b_i$. 
Let $u=\sum_i S(b_i)a_i$ be the Drinfeld element 
in $H$ (where $S$ is the antipode of $H$). Drinfeld showed 
\cite{dr1} that $u$ is invertible and 
\begin{equation}\label{u}
uxu^{-1}=S^2(u)
\end{equation}
for any $x\in H$. 
He also showed that 
\begin{equation}\label{du}
\Delta(u)=(u\ot u)(R^{21}R)^{-1}.
\end{equation}
If $H$ is also semisimple then $S^2=1$ \cite{lr}, 
hence $u$ is a central element in $H$. 

Let $H$ be a finite-dimensional semisimple Hopf algebra over $k.$ Then the 
Drinfeld double of $H,$ $D(H),$ is semisimple \cite{r} and quasitriangular
\cite{dr1} with
universal $R-$matrix $R=\sum_i 
h_i\ot h_i^*,$ where $\{h_i\}$ and $\{h_i^*\}$ are dual bases of $H$ and 
$H^*$ respectively. It is moreover 
a ribbon Hopf algebra \cite{kas} with the Drinfeld element $u,$
defined in
(\ref{u}), as the ribbon element $v.$ In particular the special grouplike
element $g=uv^{-1}$ equals $1$ in this case. Equivalently, the 
category $Rep(D(H))$ of finite-dimensional representations of $D(H)$ is a 
semisimple ribbon (i.e.
braided, rigid and balanced) category with quantum trace equal to the
ordinary trace.
Let $Irr(D(H))=\{V_i|0\le i\le m\}$ be the set of all the irreducible 
representations of $D(H)$ with $V_0=k,$ and let $C(D(H))\subseteq D(H)^*$
be the ring of characters. Clearly, 
$\{\chi _i=tr_{|V_i}|0\le i\le m\}$ forms a linear basis of $C(D(H)).$
We also let $\chi _{i^*}=S(\chi _i)$ be the character of the irreducible
representation $V_{i^*}=V_{i}^*.$

Recall that a {\em modular} category \cite{ki} is a semisimple ribbon
category with
finitely many irreducible objects $\{V_i|0\le i\le m\}$ so that the matrix
$s=(s_{ij}),$ where $s_{ij}=(tr_{|{V_i}}\ot tr_{|{V_{j^*}}})(R^{21}R)$ is
invertible. Note that $s$ is symmetric and $s_{i0}=s_{0i}=dimV_i$ for all $i.$
\begin{Lemma}\label{modular} 
Let $H$ be a finite-dimensional semisimple Hopf algebra over $k.$ Then 
$Rep(D(H))$ is a modular category. 
\end{Lemma}
\proof We only have to show that the matrix $s=(s_{ij}),$ where $s_{ij}= 
(\chi _i\ot \chi _{j^*})(R^{21}R)$ is invertible.
Indeed, by \cite{dr1}, $D(H)$ is factorizable (i.e. the map $F:D(H)^*\raro
D(H)$ given by $F(p)=(1\ot p)(R^{21}R)$ is an isomorphism of
vector spaces), and $F_{|C(D(H))}: C(D(H))\raro Z(D(H))$ is an
isomorphism of algebras, where $Z(D(H))$ is the center of $D(H).$ 
Let $B=\{\chi _{j^*}|0\le j\le m\}$ and $C=\{e_j|0\le j\le m\}$ be bases
of $C(D(H))$ and $Z(D(H))$ respectively, where $C$ is the set of central
primitive idempotents of $D(H).$ Then it is 
straightforward to check that $s=AD$ where $A$ is
the invertible matrix which represents $F$ with respect to the bases $B$
and $C,$ and 
$D=diag(dimV_i)$ is the invertible diagonal matrix with entries $dimV_i.$
Thus, $s$ is invertible.\qed
\begin{Lemma}\label{divisible}
Let ${\cal C}$ be a modular category over $k$ with irreducible 
objects $\{V_i|0\le i\le m\}$ with $V_0=1.$ Set $R=\bigoplus
_{i=0}^{m}V_i\ot V_i^*.$ Then 
$\frac{dimR}{(dimV_j)^2}$ is an algebraic integer for all $0\le j\le m.$
\end{Lemma}
\proof 
It is known that $\sum_i s_{ji}s_{ij^*}=dimR$ for all $j$ (see e.g. \cite
{ki}), hence
$\sum\frac{s_{ji}}{s_{j0}}\frac{s_{ij^*}}{s_{0j^*}}=\frac{dimR}{(dimV_j)^2}.$
We show that $\frac{s_{ij}}{s_{i0}}$ is an algebraic integer for all
$0\le i,j\le m.$ Define a map $\phi$ from $Rep({\cal C})$ to the
algebra of functions $k\{i|0\le i\le m\}$ by
$$\phi(V_j)(i)=\frac{s_{ij}}{s_{i0}}=\frac{1}{dimV_i}
(tr_{|{V_i}}\ot tr_{|{V_{j^*}}})(R^{21}R).$$ It is straightforward to check
that $\phi$ is an isomorphism of algebras (see e.g.
\cite{t}). Since multiplication by
$\phi(V_j)$ has eigenvalues $\{\frac{s_{ij}}{s_{i0}}|0\le i\le m\}$ it
follows that multiplication by $V_j$ in $Rep({\cal C})$
($V_i\mapsto V_j\ot V_i$) has the same eigenvalues (this statement is
called "Verlinde Formula" \cite{v}). But,
multiplication is represented by an integral matrix $(N_{ij}^{l})_{jl}$
where 
$V_j\ot V_i=\sum_lN_{ij}^{l}V_l.$ We thus conclude that 
$\frac{s_{ji}}{s_{j0}}$ is an algebraic integer. Since $s$ is symmetric it 
follows that $\frac{s_{ij^*}}{s_{0j^*}}$ is an algebraic integer too.\qed
\begin{Remark}\label{eigen}
{\rm We demonstrate that the map $\phi:Rep(D(H))\raro k\{i|0\le i\le m\}$
is an algebra map. Indeed,
$\phi(V_j)(i)=\frac{1}{dimV_i}\chi_i(F(\chi_{j^*}))$ where $F$ is as in
the proof of Lemma \ref{modular}. Since $F_{|C(D(H))}$ is an isomorphism
of algebras onto $Z(D(H))$ we have 
\begin{eqnarray*}
\lefteqn{\phi(V_j\ot
V_l)(i)=\frac{1}{dimV_i}\chi_i(F(\chi_{j^*}\chi_{l^*}))=
\frac{1}{dimV_i}\chi_i(F(\chi_{j^*})F(\chi_{l^*}))}\\
& & =\frac{1}{(dimV_i)^2}\chi_i(F(\chi_{j^*}))\chi_i(F(\chi_{l^*}))
=\phi(V_j)(i)\phi(V_l)(i).
\end{eqnarray*}
}
\end{Remark}
\begin{Theorem}\label{main2}
Let $H$ be a finite-dimensional semisimple Hopf algebra over $k$. If $V$ is
an irreducible representation of $D(H)$ then $dimV$ divides $dimH.$
\end{Theorem}
\proof First note that since $D(H)$ is semisimple, $D(H)=\bigoplus
_{i=0}^{m}V_i\ot V_i^*.$ Now, 
by Lemmas \ref{modular} and \ref{divisible},
$\frac{(dimH)^2}{(dimV)^2}$ is an algebraic integer. This implies that
$\frac{dimH}{dimV}$ is an algebraic integer too, hence an integer.\qed

We are ready now to prove Kaplansky's conjecture for quasitriangular
semisimple Hopf
algebras.
\begin{Theorem}\label{main3}
Let $(H,R)$ be a finite-dimensional quasitriangular semisimple Hopf
algebra. Then $H$ is of Frobenius type. 
\end{Theorem}
\proof By \cite {dr1}, the map $f:D(H)\raro H$ given by $f(ph)=(p\ot
1)(R)h$ for all $p\in H^*$ and $h\in H,$ is a surjection of Hopf algebras.
Therefore, if $V$ is an irreducible
representation of $H$ then it
is also an irreducible representation of $D(H)$ via pull back along $f,$
and the result follows by Theorem \ref{main2}.\qed
\begin{Example}\label{gral} {\rm The group algebra $H=kG$ of a
finite group $G$ is quasitriangular with $R=1\ot 1.$ 
In this case Theorem \ref{main3} is the classical theorem stating that
the dimensions 
of the irreducible representations of a finite group divide its order. 
In fact, our proof of Theorem \ref{main3} in the case of group algebras
reproduces one of the classical proofs of this well-known result.}
\end{Example}
In the following we show how Kaplansky's conjecture on prime
dimensional Hopf algebras follows
easily from Theorem \ref{main2} (compare with \cite{z}).
\begin{Theorem}\label{p}
Let $H$ be a Hopf algebra of prime dimension $p.$ Then $H=k\Z_p$ is the
group algebra of the cyclic group of order $p.$
\end{Theorem}
\proof By \cite{nz}, either $|G(H)|=1$ or $|G(H)|=p,$ and the same holds
for $H^*.$ Suppose $|G(H)|=|G(H^*)|=1.$ Then it is easy to show that 
$H$ is semisimple (see e.g.\cite{z}).
But then by Theorem
\ref{main2}, if $V$ is an irreducible representation of $D(H)$ then either
$dimV=1$ or $dimV=p.$ 
Since $G(D(H)^*)=1$ (i.e. $D(H)$ has only one 
$1-$dimensional representation) and $dimD(H)=p^{2},$ 
it follows that $p^{2}=1+ap^{2}$ for some positive integer $a$ which is 
absurd. Therefore, either $|G(H^*)|=p$ or $|G(H)|=p,$ and the result
follows.\qed
\section{Triangular Hopf algebras}
By Theorem \ref{main3}, if $(H,R)$ is a finite-dimensional triangular
(i.e. $R^{21}R=1$) semisimple Hopf
algebra then it is of Frobenius type. In fact, we can say much more in
this case. It was conjectured in \cite{g} that the minimal sub Hopf
algebra $H_R\subseteq H$ is
isomorphic to a group algebra as a Hopf algebra. In the following we show
that $H$ itself is a twisted group algebra in the sense of \cite{dr2}; 
that is, that $H$ is isomorphic as a Hopf algebra to a group algebra with a 
deformed comultiplication of the standard one by an invertible counital 
cocycle.

Let $(H,R)$ be a finite-dimensional semisimple triangular Hopf
algebra. In this case $u$ is a grouplike element by (\ref{du}).  
\begin{Lemma} $u^2=1$.
\end{Lemma}
\proof We have $(S\ot S)(R)=R$, so $u^{-1}=S(u)=
\sum S(a_i)S^2(b_i)=\sum a_iS(b_i)$. This shows that $tr(u)=tr(u^{-1})$ in
 every irreducible representation of $H$. But $u$ is central, so it 
acts as a scalar in this representation. Thus, $u=u^{-1}$. 
\qed

For some purposes it is useful to assume that the Drinfeld element $u$ 
acts as 
$1$ (see e.g. \cite{cwz}). Let us demonstrate that it is always possible 
to replace $R$ with a new $R-$matrix $\tilde R$ so that the new Drinfeld 
element $\tilde u$ equals $1.$ Indeed, for any irreducible 
representation $V$ of $H$ define the parity of this representation, 
$p(V)\in \Z_2$, by $(-1)^{p(V)}=u|_{V}$. 
Define $\tilde R\in H\ot H$ by the condition 
$\tilde R|_{V\ot W}=(-1)^{p(V)p(W)}R|_{V\ot W}.$ 
It is straightforward to verify that $\tilde R=\frac{1}{2}(1\ot 1+1\ot 
u+u\ot 1-u\ot u)R$ is a new triangular structure on $H$, with Drinfeld 
element $\tilde u=1$. 

Our main result in this section is:
\begin{Theorem}\label{main1} 
Let $(H,R)$ be a finite-dimensional semisimple triangular Hopf algebra
over $k.$ 
Then there exists a finite group $G$, an invertible element
$J\in kG\ot kG$ which satisfies 
\begin{equation}\label{cocycle}
(\Delta_0\ot 1)(J)J_{12}=(1\ot \Delta_0)(J)J_{23}, 
\,\,\,(\varepsilon_0\ot 1)(J)=(1\ot \varepsilon_0)(J)=1,
\end{equation}
(where $\Delta_0,\varepsilon_0$ are the coproduct 
and the counit of the group algebra), and 
an algebra isomorphism $\phi:kG\raro H$ such that 
\begin{equation}\label{gauge}
(\phi^{-1}\ot \phi^{-1})(\Delta(\phi(a)))=J^{-1}\Delta_0(a)J,
\end{equation}
and 
\begin{equation}\label{twistr}
(\phi^{-1}\ot \phi^{-1})(\tilde R)=(J^{21})^{-1}J.
\end{equation} 
That is, $(H,\tilde R)$ and $(kG,\gD,(J^{21})^{-1}J)$ are isomorphic as  
triangular Hopf algebras, where $\gD :kG\raro kG\ot kG$ is
determined by $\gD (g)=J^{-1}(g\ot g)J.$
\end{Theorem}
\proof Let ${\cal C}$ be the category of finite-dimensional
representations of $H$. 
This is a semisimple abelian 
category over $k$ with finitely many irreducible objects, 
which has a structure of a rigid symmetric tensor
category \cite{dm}. Here the commutativity isomorphism in ${\cal C}$ is
defined by 
the operator $\tau\tilde R:V\ot W\raro W\ot V,$ where $\tau :V\ot W\raro
W\ot V$ is the usual permutation map. Moreover, 
the categorical dimension \cite{dm} of an object $V\in {\cal C}$ is equal
to $tr|_V(\tilde u)$, so it equals to the ordinary dimension of
$V$ as a vector space (since $\tilde u=1$). In particular, all categorical 
dimensions are non-negative integers. 

In this situation we can apply the following deep 
theorem of Deligne:

\noin {\bf Theorem [De, Theorem 7.1]}
Let ${\cal C}$ be a semisimple rigid symmetric tensor category
over an algebraically closed field $k$ with finitely many irreducible
objects, in which 
categorical dimensions of objects are nonnegative integers.  
Then for a suitable finite group $G$ 
there exists an equivalence of symmetric rigid tensor categories 
$F:{\cal C}\raro Rep(G)$ (where $Rep(G)$ is the category of finite
dimensional 
$k-$representations of $G$).

So let $G,F$ be the group and the functor corresponding 
to our category ${\cal C}$. 
Let $K:Rep(G)\raro Vect$ be the forgetful
functor to the category of vector spaces. 
Since the functor $F$ preserves dimensions, and the category 
is semisimple, the functor $K\circ F$ is isomorphic 
(as an additive functor) to the forgetful functor 
$L:{\cal C}\raro Vect$. We might as well assume that $K\circ F=L$
as additive functors.  

By a standard argument we have $End(L)=H$. On the other hand, 
the group $G$ by definition acts on $K\circ F$ as a tensor functor, which 
defines an algebra homomorphism $\phi: kG\raro H$. 
It is obvious that $\phi$ is an algebra isomorphism. 

The functor $K\circ F$ has a tensor structure which preserves the 
commutativity isomorphism. This structure 
is given by a collection of 
invertible linear maps 
$$\tilde J_{VW}:(K\circ F)(V)\ot (K\circ F)(W)\raro 
(K\circ F)(V\ot W)$$ 
for irreducible $V,W$,
which can be united in an invertible element $\tilde J\in End((K\circ
F)^2)=
kG\ot kG$  
(since $(K\circ F)(V)\ot (K\circ F)(W)=L(V)\ot L(W)=L(V\ot W)
=(K\circ F)(V\ot W)$). 
The element $J=\tilde J^{-1}$ satisfies (\ref{cocycle}) and (\ref{gauge})
because $\tilde J$ is a tensor structure, 
and satisfies (\ref{twistr}) because $\tilde J$ preserves the
commutativity isomorphism. \qed
\begin{Remark}\label{dim} {\rm One should distinguish between 
the categorical dimensions of objects, 
defined in any rigid braided tensor category, and their 
quantum dimensions, defined only in a ribbon category. 
In the diagrammatic language of \cite{kas,ki} the quantum dimension 
corresponds to a loop without self-crossing, and the categorical 
dimension to a loop with one self-overcrossing. They may be different numbers 
for a particular irreducible object. For example, in the category of 
representations of a triangular semisimple Hopf algebra $(H,R)$, 
quantum dimensions (for an appropriate ribbon structure) are ordinary 
dimensions (as in Section 1), 
while categorical dimensions are $u|_V\text{dim}(V)$, where 
$u|_{V}$ is the scalar by which the Drinfeld element $u$ acts on $V,$
i.e. $1$ or $-1$ (as in Section 2).} 
\end{Remark}
\begin{Remark}\label{sv} {\rm As seen from Remark \ref{dim}, if $u\ne
1$, then the
category of representations of $(H,R)$ is equivalent to the category 
of representations of some group as a rigid tensor category 
but not as a symmetric category. This was the reason for passing 
from $R$ to $\tilde R$. It is easy to see that as a symmetric
rigid tensor category, the category of representations of $(H,R)$ is 
equivalent to the category of representations of $G$ on super-vector spaces, 
such that $u$ acts by $1$ on the even part and as $-1$ on the odd part. 
For example, if $H=k\Z_2$ with central primitive idempotents
$a$ and $b,$ and $R=a\otimes a+b\otimes a+a\otimes b
-b\otimes b$, then the category of representations is just the category 
of super-vector spaces.} 
\end{Remark}
\noin
{\bf Acknowledgments} The first author is grateful to 
Dennis Gaitsgory and Alexander Kirillov Jr. for useful discussions.
The first author was supported by an NSF postdoctoral fellowship.
The second author was supported by a Rothschild postdoctoral fellowship.
 
\end{document}